\documentclass[twoside,psfig]{article}
\usepackage{amsfonts}
\usepackage{latexsym}

\usepackage{times}
\usepackage{graphicx}
\usepackage{epsfig}
\usepackage{amsmath}

\newcommand{\addeps}[3]  
{\begin{figure}[htbp]
\begin{center}
\epsfig{file=#1.eps,width=#2in}
\caption{#3}\label{#1}
\end{center}\end{figure}}

\oddsidemargin=\evensidemargin 
\newcommand{\Z}{{\mathbb Z}}
\newcommand{\sii}{\sigma_i}
\newcommand{\sij}{\sigma_j}
\newcommand{\ti}{\tau_i}
\newcommand{\tj}{\tau_j}
\catcode`\@=11 \long\def\@makefntext#1{ \protect\noindent \hbox to
3.2pt {\hskip-.9pt
$^{{\eightrm\@thefnmark}}$\hfil}#1\hfill}       

\def\ps@myheadings{\let\@mkboth\@gobbletwo      
\def\@oddhead{\hbox{}
\rightmark\hfil\eightrm\thepage}
\def\@oddfoot{}\def\@evenhead{\eightrm\thepage\hfil
\leftmark\hbox{}}\def\@evenfoot{}
\def\sectionmark##1{}\def\subsectionmark##1{}}

\catcode`@=11                        
\def\ps@plain{\let\@mkboth\@gobbletwo
     \def\@oddhead{}\def\@oddfoot{\eightrm\hfil\thepage
     \hfil}\def\@evenhead{}\let\@evenfoot\@oddfoot}

\newcounter{sectionc}\newcounter{subsectionc}\newcounter{subsubsectionc}
\renewcommand{\section}[1] {\vspace{12pt}\addtocounter{sectionc}{1}
\setcounter{subsectionc}{0}\setcounter{subsubsectionc}{0}\noindent
    {\tenbf\thesectionc. #1}\par\vspace{5pt}}
\renewcommand{\subsection}[1] {\vspace{12pt}\addtocounter{subsectionc}{1}
    \setcounter{subsubsectionc}{0}\noindent
    {\bf\thesectionc.\thesubsectionc.
    {\kern1pt \bfit #1}}\par\vspace{5pt}}
\renewcommand{\subsubsection}[1] {\vspace{12pt}
    \addtocounter{subsubsectionc}{1}
    \noindent
    {\tenrm\thesectionc.\thesubsectionc.\thesubsubsectionc. {\kern1pt
    \it #1}}\par\vspace{5pt}}

\newcounter{appendixc}
\newcounter{subappendixc}[appendixc]
\newcounter{subsubappendixc}[subappendixc]

\renewcommand{\appendix}[1] {\vspace{12pt}  
    \refstepcounter{appendixc}      
    \setcounter{figure}{0}
    \setcounter{table}{0}
    \setcounter{lemma}{0}
    \setcounter{theorem}{0}
    \setcounter{corollary}{0}
    \setcounter{definition}{0}
    \setcounter{equation}{0}
    \renewcommand{\thefigure}{\Alph{appendixc}.\arabic{figure}}
    \renewcommand{\thetable}{\Alph{appendixc}.\arabic{table}}
    \renewcommand{\theappendixc}{\Alph{appendixc}}
    \renewcommand{\thelemma}{\Alph{appendixc}.\arabic{lemma}}
    \renewcommand{\thetheorem}{\Alph{appendixc}.\arabic{theorem}}
    \renewcommand{\thedefinition}{\Alph{appendixc}.\arabic{definition}}
    \renewcommand{\thecorollary}{\Alph{appendixc}.\arabic{corollary}}
    \renewcommand{\theequation}{\Alph{appendixc}.\arabic{equation}}
    \noindent{\tenbf Appendix \theappendixc #1}\par\vspace{5pt}}

\newcommand{\smalllineskip}{\baselineskip=10pt}

\newcommand{\copyrightheading}[1]
    {\vspace*{-2.5cm}\smalllineskip{\flushleft
    {\footnotesize Journal of Knot Theory and Its Ramifications #1}\\
    {\footnotesize \copyright\kern2pt World Scientific
         Publishing Company}\\
         }}

\def\keywords#1{{
    \centering{\begin{minipage}{4.5in}\footnotesize\baselineskip=10pt
    {\footnotesize\it Keywords}\/: #1
    \end{minipage}}\par}}

\renewenvironment{thebibliography}[1]
    {\frenchspacing
     \ninerm\baselineskip=11pt
     \begin{list}{[\arabic{enumi}]}
    {\usecounter{enumi}\setlength{\parsep}{0pt}
     \setlength{\leftmargin 13.7pt}{\rightmargin 0pt} 
     \setlength{\itemsep}{0pt} \settowidth
    {\labelwidth}{[#1]}\sloppy}}{\end{list}}

\newcounter{itemlistc}
\newcounter{romanlistc}
\newcounter{alphlistc}
\newcounter{arabiclistc}

\newcommand{\fcaption}[1]{
        \refstepcounter{figure}
        \setbox\@tempboxa = \hbox{\footnotesize Fig.~\thefigure. #1}
        \ifdim \wd\@tempboxa > 5in
           {\begin{center}
        \parbox{5in}{\footnotesize\smalllineskip Fig.~\thefigure. #1}
            \end{center}}
        \else
             {\begin{center}
             {\footnotesize Fig.~\thefigure. #1}
              \end{center}}
        \fi}

\newcommand{\tcaption}[1]{
        \refstepcounter{table}
        \setbox\@tempboxa = \hbox{\footnotesize Table~\thetable. #1}
        \ifdim \wd\@tempboxa > 5in
           {\begin{center}
        \parbox{5in}{\footnotesize\smalllineskip Table~\thetable. #1}
            \end{center}}
        \else
             {\begin{center}
             {\footnotesize Table~\thetable. #1}
              \end{center}}
        \fi}

\def\pmb#1{\setbox0=\hbox{#1}
    \kern-.025em\copy0\kern-\wd0
    \kern.05em\copy0\kern-\wd0
    \kern-.025em\raise.0433em\box0}

\def\fnt#1#2{\footnotetext{\kern-.3em
    {$^{\mbox{\scriptsize #1}}$}{#2}}}

\font\tenrm=cmr10  \font\tenbf=cmbx10
\font\bfit=cmbxti10 at 10pt \font\ninerm=cmr9 
 \font\eightrm=cmr8

\newtheorem{thm}{Theorem}   

\newtheorem{lemma}{Lemma}

\newtheorem{remark}{Remark}

\def\@begintheorem#1#2{\trivlist    
    \item[\hskip\labelsep{\bf #1\ #2.}]}
\def\@opargbegintheorem#1#2#3{\trivlist
    \item[\hskip\labelsep{\bf #1\ #2\ (#3).}]}

\newenvironment{proof}{\begin{trivlist}
    \item[\noindent]{\it Proof.}}{\end{trivlist}}

        {\setcounter{itemlistc}{0}      
     \begin{list}{$\bullet$}        
    {\usecounter{itemlistc}         
     \leftmargin10pt           
     \setlength{\parsep}{0pt}
     \setlength{\itemsep}{0pt}     
    }}{\end{list}}

    {\setcounter{romanlistc}{0}     
     \begin{list}{$($\roman{romanlistc}$)$} 
    {\usecounter{romanlistc}        
     \leftmargin18pt 
     \setlength{\parsep}{0pt}
     \setlength{\itemsep}{0pt}  
     \settowidth{\labelwidth}{#1}
    }}{\end{list}}

    {\setcounter{enumii}{0}         
     \begin{list}{$($\alph{enumii}$)$}  
    {\usecounter{enumii}            
     \leftmargin18pt        
     \setlength{\parsep}{0pt}
     \setlength{\itemsep}{0pt}  
     \settowidth{\labelwidth}{#1}
    }}{\end{list}}

\def\qed{\hbox{${\vcenter{\vbox{            
   \hrule height 0.4pt\hbox{\vrule width 0.4pt height 6pt
   \kern5pt\vrule width 0.4pt}\hrule height 0.4pt}}}$}}

\def\theequation{\thesectionc.\arabic{equation}}  
\begin{document}
\setlength{\textheight}{7.7truein} \centerline{\bf Classifying links under fused isotopy} \vspace*{0.37truein} \centerline{\footnotesize ANDREW
FISH and EBRU KEYMAN} \baselineskip=12pt
\vspace*{5pt}\centerline{\footnotesize\it School of Computing,
Mathematical and Information Sciences, University of Brighton, UK
} \baselineskip=10pt \centerline{\footnotesize\it
Andrew.Fish@brighton.ac.uk} \vspace*{10pt} \baselineskip=12pt
\centerline{\footnotesize\it Department of Mathematics, Middle
East Technical University, Ankara, Turkey} \baselineskip=10pt
\centerline{\footnotesize\it ekeyman@metu.edu.tr}

\vspace*{0.225truein}
\begin{abstract}
All knots are fused isotopic to the unknot using a process known
as virtualization. We extend and adapt this process to show that,
up to fused isotopy, classical links are classified by their
linking numbers.
\end{abstract}

\section{Introduction}

\begin{figure}[h!]
\centerline{\begin{picture}(260,46) 
\put(0,0){\includegraphics[scale=1]{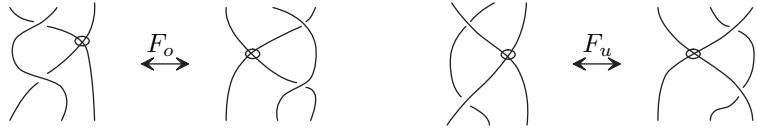}}
\put(52,28){$F_o$}\put(217,28){$F_u$}\end{picture}} \caption{The
forbidden moves}\label{forbidden}
\end{figure}

Classical braids and links have been generalized to the virtual
category~\cite{Kf} -- adding virtual crossings and extending
isotopy to allow the virtual analogues of the classical
Reidemeister moves. The forbidden moves $F_o$ and $F_u$, shown in
Figure~\ref{forbidden}, are not allowable under virtual isotopy.
Extending virtual isotopy in the virtual braid group $VB_n$ to
allow the $F_o$ move gives rise to the welded braid group $WB_n$,
which has been shown to be isomorphic to $PC_n$, the group of
automorphisms of the free group on $n$ elements of
permutation-conjugacy type~\cite{FRR}. Allowing both of the
forbidden moves $F_o$ and $F_u$ gives rise to fused
isotopy~\cite{Kf}.
%
%
%
%
%
%
That is, two virtual links $L_1$ and $L_2$ are called \emph{fused
isotopic} if $L_2$ can be obtained from $L_1$ by a finite sequence
of Reidemeister moves, virtual moves and $F_o$, $F_u$ moves.

Let $\mathcal{K}$ be the space of classical links embedded in
$S^3$ and let $\mathcal{VK}$ be the space of virtual links.
Kauffman~\cite{Kf}, and independently
Goussarov-Polyak-Viro~\cite{GPV}, have shown that $\mathcal{K}$
embeds into $\mathcal{VK}$. Let $f$ denote the natural inclusion
of $\mathcal{VK}$ into the space of fused links $\mathcal{FK}$.
Then we have
$\mathcal{K}\stackrel{i}{\hookrightarrow}\mathcal{VK}\stackrel{f}{\rightarrow}\mathcal{FK}$,
and when we refer to a classical link (under fused isotopy) we
mean $f\circ i(L)$, the image of a link $L \in \mathcal{K}$ in the
space $\mathcal{FK}$.

In~\cite{Kn}, Kanenobu showed that all knots are fused isotopic to
the unknot. He showed that all of the classical crossings of a
virtual knot can be \emph{virtualized}; that is every classical
crossing can be changed into a virtual crossing by applying a
sequence of fused isotopy moves. However, crossings between
different components of a link cannot be virtualized using the
same methods. The following theorem from~\cite{Kn} provides us with allowable
moves under fused isotopy which were used in the virtualization
procedure.

\begin{figure}[h!]
\centerline{\begin{picture}(250,108) 
\put(0,0){\includegraphics[scale=1]{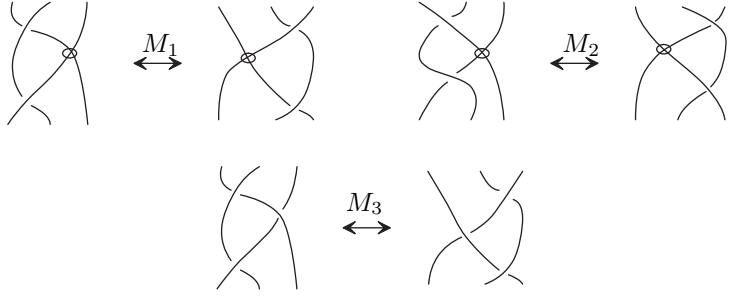}}
\put(51,90){$M_1$}\put(210,90){$M_2$}\put(128,30){$M_3$}\end{picture}}
\caption{Allowable moves in fused isotopy}\label{fused moves}
\end{figure}

\begin{thm}The moves $M_1,M_2$ and $M_3$, shown in Figure~\ref{fused moves}, can be realised by fused isotopy. \quad $\square$
\label{thm:MMoves}
\end{thm}

In~\cite{FK} the authors showed that the Jones polynomial for
welded and fused links is well-defined in a quotient of
${\ensuremath {\mathbb{Z}}}[A,A^{-1}]$ and observed that this
polynomial depends only upon the linking number for links with two
components. Inspired by this, we show that classical links, under
fused isotopy, can be determined by the linking number of their
components.

\begin{thm} A classical link $L$ with $n$--components is completely determined by the linking numbers
of each pair of components under fused isotopy. \label{MainThm}
\end{thm}

The strategy that we use to prove Theorem~\ref{MainThm} is to
write $L$ as the closure of a braid $\alpha$ on $m$ strands (where
$m \geq n$) and then to transform $\alpha$ into a pure braid
$\beta$ on $n$ strands whose closure is also $L$. We show that
$\beta$ depends only on the linking numbers of the components of
$L$. This means that any classical link with the same linking
numbers as $L$ can be obtained as the closure of $\beta$. We need
some preliminaries before we proceed with the proof.

\section{Preliminaries}
Recall that an element of the pure braid group $P_n$ is an
$n$--strand braid where the permutation induced by the strings is
the identity. $P_n$ has a presentation with generators $A_{i,j}$
with $1\leq i< j\leq n$ where
$$A_{i,j}=
\sigma_{j-1}\sigma_{j-2}\ldots\sigma_i^2\ldots\sigma_{j-2}^{-1}\sigma_{j-1}^{-1}
=
\sigma_i^{-1}\sigma_{i+1}^{-1}\ldots\sigma_{j-1}^2\ldots\sigma_{i+1}\sii.$$

\begin{figure}[ht]
\centerline{\begin{picture}(100,65)
\put(0,0){\includegraphics[scale=1.1]{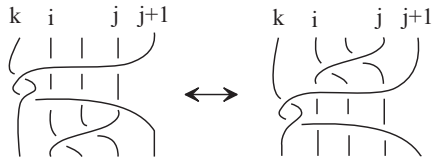}}\end{picture}}
\caption{$A_{k,j+1}B_{i,j}=B_{i,j}A_{k,j+1}$}\label{AB} 
\end{figure}

Let $U_k$ be the subgroup of $P_n$ generated by $\{A_{i,k}:\ 1\leq
i<k \}$. Then every element of $P_n$ can be written in the unique
normal form $x_2x_3\ldots x_n$, where $x_k\in U_k$ (see \cite{Bi}
for details).
Define $B_{i,j}:= \sigma_{j-1}\ldots\sigma_{i+1}\sii$ for $i<j$
and $B_{i,i}:=1$. Then by definition
$A_{i,j+1}=B_{i,j}^{-1}A_{j,j+1}B_{i,j}$, and we can see from
Figure~\ref{AB} that for $k<i<j,\ B_{i,j}$ commutes with
$A_{k,j+1}$ in $B_n$.

The virtual braid group on $n$--strands, $VB_n$, can be defined by
adding extra generators $\tau_i$, for $1 \leq i \leq n$,
corresponding to the virtual crossings, and relations
corresponding to the virtual isotopy moves~\cite{Km, Kf}. By
adding the relations $\sii^{-1}\tj\sii = \sij\ti\sij^{-1}$
%
with $\left|i-j\right|=1$, to the virtual braid group, we obtain
the fused braid group, $FB_n$. If $j=i+1$ the relation corresponds
to the $F_u$ move, and if $i=j+1$ it corresponds to the $F_o$
move.
The explicit realization of the moves $M_1, M_2$ and $M_3$ using
$F_o$ and $F_u$ moves is shown in~\cite{Kn}, and this gives rise
to the following consequences in $FB_n$:
\begin{equation}
\begin{array}{lrcl}
M_1:&\sii\tj\sii&=&\sij\ti\sij\\
M_2:&\sii^{-1}\tj\sii^{-1}&=&\sij^{-1}\ti\sij^{-1}\\
M_3:&\sii\sij^{-1}\sii & = & \sij\sii^{-1}\sij ,
\end{array}
\end{equation}
\noindent where $\left|i-j\right|=1$.

The following lemmas are used in the proof of
Theorem~\ref{MainThm}, and the indices have been chosen to match
the usage in the proof. Let $\sim$ denote the equivalence class
generated by fused isotopy and let $\widetilde{U}_k = U_k/\!
\sim$.
%
\begin{lemma}\label{central} In $FP_n$ we have:
\begin{equation}\label{com}
A_{j,j+1}A_{i,j+1}A_{j,j+1}^{-1}=A_{i,j+1}\ \mathrm{where}\ 1\leq
i<j+1\leq n.
\end{equation}
In other words, $A_{j,j+1} = \sij^2$ is in the centre of $\widetilde{U}_{j+1}$.
\end{lemma}

\begin{proof}
Using the relations
$\sij\sigma_{j-1}^{-1}\sij=\sigma_{j-1}\sij^{-1}\sigma_{j-1}$
corresponding to an $M_3$ move, and
$\sij\sigma_{j-1}\sij^{-1}=\sigma_{j-1}^{-1}\sigma_{j}\sigma_{j-1}$
corresponding to an $R_3$ move, we obtain

\begin{equation}\label{pfa}
\begin{array}{cl}
&\sigma_{j}^{2}\sigma_{j-1}^{-1}\sij^2\sigma_{j-1}\sigma_{j}^{-2}\,
=\, \sij\left( \sij\sigma_{j-1}^{-1}\sij\right)
\left( \sij\sigma_{j-1}\sigma_{j}^{-1}\right)\sij^{-1}\\
=&\sij\left( \sigma_{j-1}\sij^{-1}\sigma_{j-1}\right)
\left(\sigma_{j-1}^{-1}\sigma_{j}\sigma_{j-1}\right)\sij^{-1}\,
=\, \sij\sigma_{j-1}^2\sij^{-1} .
\end{array}
\end{equation}

Using Equation~(\ref{pfa}), we can compute
\begin{equation*}
\begin{array}{cll}
&A_{j,j+1}A_{i,j+1}A_{j,j+1}^{-1} \\ 
=&\sigma_{j}^{2}B_{i,j}^{-1}A_{j,j+1}B_{i,j}\sigma_{j}^{-2}& \\ 
%
=&B_{i,j-1}^{-1}\sij^2\sigma_{j-1}^{-1}\sij^2\sigma_{j-1}\sigma_{j}^{-2}B_{i,j-1}&\mathrm{\ by\ commutation\ in\ }B_n\\
%
=&B_{i,j-1}^{-1}\sij\sigma_{j-1}^2\sij^{-1}B_{i,j-1}&\mathrm{\ by\ Equation\ (\ref{pfa})}\\
=&\sij B_{i,j-1}^{-1}\sigma_{j-1}^2B_{i,j-1}\sij^{-1}&\mathrm{\ by\ commutation\ in\ }B_n\\
=&\sij A_{i,j}\sij^{-1}&\\ 
 =&A_{i,j+1} .&
\square
\end{array}
\end{equation*}

\end{proof}

\begin{lemma}\label{Uk commutes} For every $1\leq j+1 <n$, the subgroup $\widetilde{U}_{j+1}$ of $FP_n$ is commutative.
\end{lemma}

\begin{proof} Assume without loss of generality that $k<i$. Then
\begin{equation*}
\begin{array}{rll}
A_{k,j+1}A_{i,j+1}= &A_{k,j+1}B_{i,j}^{-1}A_{j,j+1}B_{i,j}& \\ 
= &B_{i,j}^{-1} A_{k,j+1}A_{j,j+1} B_{i,j}& \mathrm{\ by\ commutation\ in\ }B_n \\
= & B_{i,j}^{-1} A_{j,j+1}A_{k,j+1} B_{i,j}& \mathrm{\ by\ Lemma\ \ref{central}}\\
=& B_{i,j}^{-1}A_{j,j+1}B_{i,j}A_{k,j+1}& \mathrm{\ by\ commutation\ in\ }B_n\\
 =& A_{i,j+1}A_{k,j+1}.&\square
\end{array}
\end{equation*}

\end{proof}



\begin{lemma} \label{lemma:tau}
In $FB_n$ we have:
\begin{equation}\label{tau}
A_{i,j+1}\tau_j =\tau_jA_{i,j}  ~~~~~ \mathrm{\ for\ } 1\leq i\leq
j-1.
\end{equation}
\end{lemma}
\begin{proof}
Using the relations $\sij \sigma_{j-1} \tau_{j}=\tau_{j-1} \sij
\sigma_{j-1}$ corresponding to an $F_0$ move, and
$\sigma_{j-1}^{-1} \sij \tau_{j-1} = \tau_j\sigma_{j-1}\sij^{-1}$
corresponding to an $M_1$ move, we obtain
\begin{equation*}
\begin{array}{rll}
A_{i,j+1}\tau_j = & B_{i,j}^{-1}\sij^{2}B_{i,j}\tau_j & \\ 
=& B_{i,j-1}^{-1} \sigma_{j-1}^{-1}\sij \sij \sigma_{j-1} B_{i,j-1} \tau_{j}& \\
=& B_{i,j-1}^{-1} \sigma_{j-1}^{-1}\sij \sij \sigma_{j-1} \tau_{j} B_{i,j-1} & \mathrm{\ by\ commutation\ in\ }VB_n \\
=& B_{i,j-1}^{-1} \sigma_{j-1}^{-1} \sij \tau_{j-1} \sij \sigma_{j-1} B_{i,j-1}& \mathrm{\ by\ an\ }F_o\mathrm{\ move}\\
=& B_{i,j-1}^{-1} \tau_j\sigma_{j-1}\sij^{-1} \sij\sigma_{j-1}B_{i,j-1}& \mathrm{\ by\ an\ }M_1\mathrm{\ move} \\
=& \tau_j B_{i,j-1}^{-1}\sigma_{j-1}^2B_{i,j-1}&\mathrm{\ by\ commutation\ in\ }VB_n \\
= &\tau_jA_{i,j}. & \square
\end{array}
\end{equation*}
\end{proof}

\section{Proof of Theorem~\ref{MainThm}}

Let $L$ be a classical link with $n$--components. By Alexander's
Theorem, there exists $\alpha \in B_m$ with $m\geq n$, such that
the closure of $\alpha$ is $L$. Chow~\cite{C} (or see page 22 of~\cite{Bi})
shows that every $\alpha$ can be written in the form
$\alpha=x_2B_{k_2,2}\ldots x_mB_{k_m,m}$ where $x_i\in U_i \leq
P_m$ and $1\leq k_i\leq i$. Let $\hat\alpha$ denote the closure of
$\alpha$. If $m>n$ then we will construct $\beta\in B_n$ such that
$\hat\beta$ is fused isotopic to $\hat\alpha$.

\begin{figure}[h]
\centerline{\begin{picture}(40,87)
\put(0,0){\includegraphics[scale=1]{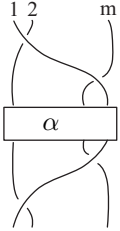}}
\put(15,37){$\alpha$}\end{picture}} \caption{$B_{1,m}^{-1}\alpha B_{1,m}$}\label{alfa'}
\end{figure}

If $B_{k_i,i}=1$ for all $i=2,\ldots ,m$ then $\alpha$ is a pure
braid, which means that $m$ must be equal to $n$. So let us assume
that $B_{k_s,s}\neq 1$, for some $s$, and that if $i>s$ then
$B_{k_i,i} = 1$. This means that the permutation induced by
$\alpha$ is the identity on the strands $s+1, \ldots, m$.
Therefore, each of these strands forms a separate component of the
link $\hat\alpha$. Now, conjugating $\alpha$ with $B_{1,m}$ gives
$B_{1,m}^{-1}\alpha B_{1,m}$, and as shown in Figure~\ref{alfa'},
the $(m-1)$--strand of the original braid $\alpha$ becomes the
$m$--strand of the new braid.

\begin{figure}[h!]
\centerline{\begin{picture}(40,74)
\put(0,0){\includegraphics[scale=1.1]{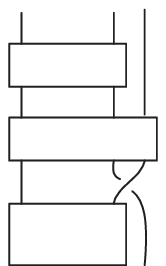}}
\put(15,62){$W$}\put(15,38){$y_m$}\put(0,8){$B_{t_m,m-1}$}\end{picture}}
\caption{The braid $\alpha'$}\label{alfa}
\end{figure}

Thus if we conjugate $\alpha$ with $B_{1,m}\ (m-s)$ times, we get
$\alpha'=B_{1,m}^{s-m}\alpha B_{1,m}^{m-s}$ and the $s$--strand of
$\alpha$ becomes the $m$--strand of $\alpha'$. Since $\alpha'$ is
just a conjugate of $\alpha$ their closures are isotopic. Now
write $\alpha'=y_2B_{t_2,2} \ldots y_m B_{t_m,m}$ with $y_i \in
\widetilde{U}_i$. Then $B_{t_m,m}\neq 1$ and so by definition,
$B_{t_{m},m} = \sigma_{m-1} B_{t_{m},m-1}$. A picture of $\alpha'$
is shown in Figure~\ref{alfa}, where $W = y_2B_{t_2,2} \ldots
y_{m-1}B_{t_{m-1},m-1}$.

Since $\widetilde{U}_m$ is commutative (by Lemma~\ref{Uk
commutes}), we can write
$$y_m=\underbrace{A_{1,m}^{r_1}\ldots A_{m-2,m}^{r_{m-2}}}_{v_{m}}
\underbrace{A_{m-1,m}^{r_{m-1}}}_{\sigma_{m-1}^{2r_{m-1}}} $$ for
some $r_1, \ldots, r_{m-1}$. By definition, $A_{m-1,m}^{r_{m-1}}
=\sigma_{m-1}^{2r_{m-1}}$, and since
$B_{t_m,m}=\sigma_{m-1}B_{t_m,m-1}$, we obtain $y_mB_{t_m,m}\
=v_{m}\sigma_{m-1}^{2r_{m-1}+1}B_{t_m,m-1}$, where $v_m =
A_{1,m}^{r_1}\ldots A_{m-2,m}^{r_{m-2}}$.

Since $W$ does not involve the $m$--strand and $y_m$ is a pure
braid, Figure~\ref{alfa} shows that the $m$--strand and the other
strand that is involved in the last occurrence (and hence in all
of the previous occurrences) of $\sigma_{m-1}$ in $\alpha'$ belong
to the same component of $L = \hat{\alpha'}$. Therefore, following
the strategy in \cite{Kn}, we can virtualize all of the
$2r_{m-1}+1$ crossings in $\hat{\alpha'}$ which correspond to
$\sigma_{m-1}^{2r_{m-1}+1}$ in $\alpha'$. In doing so we have not
changed the fused isotopy class of $L$ but we have obtained $L$ as
the closure of $\alpha_1=Wu_{m}\tau_{m-1}^{2r_{m-1}+1}B_{t_m,m-1}
= Wu_{m}\tau_{m-1}B_{t_m,m-1}$. By Lemma~\ref{lemma:tau},
we obtain $$\alpha_1=W\tau_{m-1}v_{m-1}B_{t_m,m-1}$$

\noindent where $v_{m-1} = A_{1,m-1}^{r_1}\ldots A_{m-2,m-1}^{r_{m-2}}$.

\begin{figure}[h!]
\centerline{\begin{picture}(40,74)
\put(0,0){\includegraphics[scale=1.1]{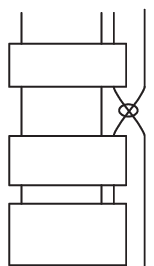}}
\put(11,59){$W$}\put(11,33){$v_{m-1}$}\put(0,8){$B_{t_m,m-1}$}
\end{picture}} \caption{The braid
$\alpha_1$}\label{newalfa}
\end{figure}

Figure~\ref{newalfa} shows that there is only one crossing
involving the $m$--strand in the braid $\alpha_1$. This is the
occurrence of $\tau_{m-1}$. In $\hat\alpha_1$, we can get rid of
the virtual crossing corresponding to $\tau_{m-1}$ with a virtual
move (of type I). We have obtained a new link diagram
$\hat\alpha_2$ where $\alpha_2=Wv_{m-1}B_{t_m,m-1}$ has $m-1$
strands
and $\hat\alpha_2$ is fused isotopic to $L$.

If we continue with this process, eventually we will get a braid $\beta$ in $B_n$ whose closure
is fused isotopic to $L$. Note that each strand of
$\beta$ corresponds to a different component of $L$ and therefore
$\beta$ must be a pure braid. For $i<j$, define the group
homomorphism $\delta_{i,j}:\ PB_n \to {\ensuremath {\mathbb{Z}}}$
by
$$
\delta_{i,j}( A_{s,t}) = \begin{cases}
1 & \text{if } s=i \text{ and } t=j \\
0 & \text{otherwise}.
\end{cases}
$$
Since $\beta$ is a pure braid
it is easy to see that $\delta_{i,j}(\beta)=lk(\ell_i,\ell_j)$
where $\ell_i$ and $\ell_j$ are the corresponding components of
$\hat{\beta}$.

This proves that any classical link $L$ with $n$--components can
be obtained as the closure of a pure braid $\beta = x_2\ldots x_n$
where each $x_k$ can be written in the form
$x_k=A_{1,k}^{\delta_{1,k}}\ldots A_{k-1,k}^{\delta_{k-1,k}}$.
This shows that $\beta$ depends only on the linking number of the
components. $\square$


\begin{remark}
We do not believe that the theorem will generalize to
non-classical links where there are virtual crossings between
different components. For example, let $U_2$ be the trivial link
with two components and let $L=\hat\alpha$ where
$\alpha=\sigma_1\tau_1\sigma_1^{-1}\tau_1\in FB_2$. Then both
 of these links have linking number 0 but we conjecture that they are not fused isotopic
 (although currently there are  no known invariants to distinguish them).
\end{remark}
%
%
%


\end{document}